\documentclass[12pt,oneside]{amsart}
\usepackage{amssymb,graphicx,amsfonts,amsthm,amsmath,latexsym}

\theoremstyle{plain}
\newtheorem{theorem}{Theorem}[section]
\newtheorem{lemma}[theorem]{Lemma}
\newtheorem{corollary}[theorem]{Corollary}
\newtheorem{proposition}[theorem]{Proposition}
\theoremstyle{definition}
\newtheorem{definition}[theorem]{Definition}

\newtheorem{example}[theorem]{Example}
\theoremstyle{remark}
\newtheorem{remark}[theorem]{Remark}

\newcommand{\R}{{\mathbb R}}
\def\K{{\bf K}}
\def\C{{\bf C}}
\def\R{{\bf R}}
\def\M{{\bf M}}
\def\H{{\bf H}}
\def\B{{\bf B}}
\def\Q{{\bf Q}}
\def\F{{\bf F}}
\def\Z{{\bf Z}}
\def\Gal{{\textrm{Gal}}}

      \makeatletter
      \def\@setcopyright{}
      \def\serieslogo@{}
      \makeatother

\begin{document}
\author[Ozman]{Ekin Ozman}
\address{Department of Mathematics,
University of Wisconsin-Madison,
480 Lincoln Drive,
Madison, WI 53706}\email{ozman@math.wisc.edu}
\title[]{Local Points on Quadratic Twists of $X_0(N)$}


\begin{abstract}

\noindent Let $X^d(N)$ be the quadratic twist of the modular curve $X_0(N)$ through the Atkin-Lehner involution $w_N$ and a quadratic extension $\Q(\sqrt{d})/\Q$.  The points of $X^d(N)(\Q)$ are precisely the $\Q(\sqrt{d})$-rational points of $X_0(N)$ that are fixed by $\sigma \circ w_N$, where 
$\sigma$ is the generator of $\Gal(\Q(\sqrt{d})/\Q)$.  
Ellenberg~\cite{Ell} asked the following question:
\begin{center}
{\it For which $d$ and $N$ does $X^d(N)$ have rational points over \\ every completion of $\Q$?}
\end{center}
Given $(N,d,p)$ we give
necessary and sufficient conditions for the existence of a $\Q_p$-rational point on 
$X^d(N)$, whenever $p$ is not simultaneously ramified in $\Q(\sqrt{d})$ and $\Q(\sqrt{-N})$, answering Ellenberg's question for all odd primes $p$ when $(N,d)=1$.
The main theorem yields a population of curves which have local points everywhere but no points over $\Q$; in several cases we show that this obstruction to the Hasse Principle is explained by the Brauer-Manin obstruction.
\end{abstract}  
\maketitle
   \section{Introduction}

Let $N=p_1 \cdots p_r$ be a positive, square-free integer. The modular curve $Y_0(N)$ is a moduli space of tuples $(E,C)$, where $E$ is an elliptic curve and $C$ is a cyclic subgroup of order $N$ in $E[N]$. Equivalently, any point of $Y_0(N)$ corresponds to $(E,\phi)$ where $\phi$ is a cyclic $N$-isogeny of $E$. A projective smooth curve $X_0(N)$ is obtained by adding $2^r$ cusps to $Y_0(N)$. 

The Atkin-Lehner involution $w_N$ of $Y_0(N)$ sends $(E,C)$ to the pair $(E/C, E[N]/C)$. Equivalently, in terms of isogenies, $w_N:  (E,\phi) \mapsto (E',\hat{\phi})$ where $\phi: E \rightarrow E'$ and $\hat{\phi}$ is the dual isogeny. The action of the rational map $w_N$ extends to $X_0(N)$ such that it permutes the cusps. 

A celebrated theorem of Mazur \cite{Mazur2} and its extensions by Kenku and Momose, give much more information about the rational points of $X_0(N)$.

\noindent
\textbf{Theorem: [Mazur, Kenku-Momose]} $X_0(N)(\Q)$ consists of only cusps when $N>163$.

This result is proved by Mazur for prime levels of $N$ and generalized to square-free integers by Kenku and Momose.
   
 Let $\K:=Q(\sqrt{d})$, $\sigma$ be the generator of $\Gal(\K/\Q)$ and $N$ be a square-free integer. The twist $X^d(N)$ of $X_0(N)$ is constructed by etale descent from $X_0(N)/\K$. It is a smooth proper curve over $\Q$,  isomorphic to $X_0(N)$ over $\K$ but not over $\Q$. The action of $\sigma$ is `twisted' on $X^d(N)$, meaning that $\Q$-rational points of $X^d(N)$ are naturally identified with the $\K$-rational points of $X_0(N)$ that are fixed by $\sigma \circ w_N$. We are interested in those points. However, the existence of rational points in this case is not immediate, as it is for $X_0(N)$. Since cusps are interchanged by $w_N$, they are not rational anymore.
 
 Like $X_0(N)$, the twisted curve $X^d(N)$ is a moduli space.  Rational points on $X^d(N)$ parametrize a special class of elliptic curves, called \emph{quadratic $\Q$-curves of degree $N$}. A $\Q$-curve is an elliptic curve which is isogenous to all of its Galois conjugates, and appears in many interesting questions, such as questions about `twisted' Fermat equations.  More details about these results and in general about $\Q$-curves, as well as related questions can be found in Ellenberg's survey article \cite{Ell}. 
 
 In general, there is no algorithm to follow when trying to show existence or non-existence of a rational point on a variety. One of the first things to check is the existence of adelic points. If a curve over $\Q$ fails to have a $\Q_p$-point for some $p$, then there is no rational point on that curve. This gives rise to the main question of the paper, which was originally stated as Problem A by Ellenberg in \cite{Ell}:
 
 \textbf{Question [Ellenberg]:} \emph{For which $d$ and $N$ does $X^d(N)$ have rational points over every completion of $\Q$?}
   

   

The main theorem of this paper gives an answer to this question under the assumption that no prime is ramified in both $\K$ and $\Q(\sqrt{-N})$.

\begin{theorem} \label{thm:ALL}
Let $p$ be a prime, $N$ square-free integer, $\K$ a quadratic field. Then 
\begin{enumerate}
\label{mainthm:real} \item \emph{$X^d(N)(\Q_p) \neq \emptyset$ for all $p$ that split in $\K$ and for $p=\infty$.}(Proposition \ref{propreal})
\label{mainthm:inertnotdividing} \item \emph{$X^d(N)(\Q_p) \neq \emptyset$ if $p$ is inert in $\K$ and not dividing $N$.}(Theorem \ref{thm:inertnotdividing})
\label{mainthm:IDB} \item \emph{For all odd $p$ inert in $\K$ and dividing $N$, $X^d(N)(\Q_p) \neq \emptyset$ if and only if}
\begin{enumerate}
	\item \emph{$N=p\Pi q_i$ where $p \equiv 3$ modulo $4$ and $q_i \equiv 1$ modulo $4$ for all $i$ or} 
	\item \emph{$N=2p\Pi q_i$ where $p \equiv 3$ modulo $4$ and $q_i \equiv 1$ modulo $4$ for all $i$.}(Theorem \ref{thm:IDB})
\end{enumerate}
\label{mainthm:inerteven}\item \emph{If $2$ is inert in $\K$ and divides $N$, $X^d(N)(\Q_p) \neq \emptyset$ if and only if $N=2\Pi q_i$ where $q_i \equiv 1$ modulo $4$ for all $i$.}(Theorem \ref{thm:inerteven})
\label{mainthm:ram}\item \emph{For all  $p$ ramified in $\K$ and not dividing $N$, $X^d(N)(\Q_p) \neq \emptyset$ if and only if $p$ is in the set $S$ defined in Proposition \ref{propmain}.} (Theorem \ref{thm:ram})
\end{enumerate}
\end{theorem}

\textbf{Organization of the paper:} In Section \ref{sec:prev} we give an overview of the previous results and where our result fits in the general scheme. We prove under the assumption that no prime is ramified in both $\K$ and $\Q(\sqrt{-N})$, the previous results on this problem by Clark, Gonzalez, Quer, and Shih are implied by Theorem \ref{thm:ALL}. We also answer a question of Clark raised in \cite{Clark}.



We draw on a number of different techniques to handle the cases in Theorem \ref{thm:ALL}.
In Section \ref{sec:real}, we deal with the case $p=\infty$.  In Section \ref{sec:inert}, we handle the case when $p$ is inert in $\K$.   For instance, in this section we use Hensel's Lemma when $p\mid N$.  On the other hand, in Section \ref{sec:ram}, we construct a $\Q_p$-point using the theory of CM elliptic curves.

In Section \ref{sec:exam}, we give some examples that illustrate Theorem \ref{thm:ALL} and compare our result with previous results.  In Section \ref{sec:further}, we give ideas about further directions, give examples of genus $2$ curves that violate Hasse principle, and show that these violations are explained by the Brauer-Manin obstruction.

\textbf{Acknowledgements:} This paper is part of my Ph.D. thesis, to be submitted in Spring 2010.  I would like to thank my advisor, Jordan Ellenberg, for suggesting the main problem of the paper, and for advice about technical matters and about the presentation.  I would also like to thank Pete Clark for several conversations on the topic, Robert Rhoades and Steven Galbraith for useful comments on the paper.

   \section{Relation with Previous Work}\label{sec:prev}
   \label{prev}
   
    
   In the case of conics 
there is a rational point if and only if there is a local point for every completion of $\Q$. Moreover if a conic has a rational point then it has many others, since it is isomorphic to $\mathbb{P}^1$. Hence for the case of conics i.e. when genus of $X^d(N)$ is zero, we have a complete answer to this question due to the work of Shih \cite{Shih}, Gonzalez and Quer given in \cite{Quer} based on the earlier work of Hasegawa \cite{Has}. The proof is based on a special parametrization of the $j$-invariants of these curves and some Hilbert symbol computations. Hence Gonzalez and Quer give the following complete list in the case of genus 0. 
   
   \begin{theorem}[Gonzalez,Quer,Shih]\label{thm:conic}
  Using the notation above:
   \begin{itemize}
   \item When $N=2,3,7$: $X^d(N)(\Q)$ is infinite for any quadratic field $\Q(\sqrt{d})$.
   \item When $N=5$    : $X^d(N)(\Q)$ is infinite if and only if $d$ is of the form $m$ or $5m$ where $m$ is a square-free integer each of whose prime divisors are quadratic residues modulo $5$.
   \item When $N=6$    : $X^d(N)(\Q)$ is infinite if and only if $d$ is of the form $m$ or $6m$ where $m$ is a square-free integer such that $2$ is a quadratic residue modulo each prime divisor of $m$.
   \item When $N=10$   : $X^d(N)(\Q)$ is infinite if and only if $d$ is of the form $m$ or $10m$ where $m$ is a square-free integer each of whose prime divisors are quadratic residues modulo $5$.
   \item When $N=13$   : $X^d(N)(\Q)$ is infinite if and only if $d$ is of the form $m$ or $13m$ where $m$ is a square-free integer each of whose prime divisors are quadratic residues modulo $13$.
   \end{itemize}
   \end{theorem}


Note that since $X^d(N)$ and $X_0(N)$ are isomorphic over $\K$ they are geometrically the same in particular they have the same genus. Therefore the cases that we know the answer completely corresponds to the values $N=2,3,5,6,7,10$ and $13$.

For $N=2,3,7$, since the class number of $\Z[\sqrt{-N}]$ is $1$, any $w_N$-fixed point of $X_0(N)$ is defined over $\Q$, hence gives a point in $X^d(N)(\Q)$ for any $d$. This is another way of saying the first part of Theorem \ref{thm:conic}. Now we'll derive the other parts of Theorem \ref{q} using Theorem \ref{thm:ALL} for relatively prime $N$ and $d$.

\begin{corollary} Let $N,d$ be square-free integers such that there is no prime $p$ that is ramified in both $\Q(\sqrt{-N})$ and $\Q(\sqrt{d})$. Then Theorem \ref{thm:conic} can be derived from Theorem \ref{thm:ALL}.
\end{corollary}

\begin{proof} 
Since we are dealing with the conics having a $\Q$-rational point is equivalent to have $\Q_p$ points for every prime $p$. By Proposition $\ref{propreal}$ $X^d(N)(\R)\neq \emptyset$ for any $N$ and $d$ hence we only need to check the finite primes. Let $d=\pm \Pi_ip_i$ be the prime decomposition of $d$.
\begin{itemize}
	\item $N=5$: By Theorem $\ref{thm:ALL}$ part 5, $X^d(5)(\Q_{p_i}) \neq \emptyset$ if and only if there is a prime of $\Q(j(\sqrt{-5}))$ lying over $p$ with inertia degree $1$. Note that class number of $\Z[\sqrt{-5}]$ is $2$ and it is the maximal order of $\M:=\Q(\sqrt{-5})$. Hilbert class field of $\M$ is $\Q(\sqrt{5},i)$ hence $\Q(j(\sqrt{-5}))$ is $\Q(\sqrt{5})$, since $j(\sqrt{-5})$ is real. Therefore $X^d(5)(\Q_{p_i}) \neq \emptyset$ if and only if $ ( \frac{5}{p_i}  )=1$.

	Since $ ( \frac{5}{p_i}  )= ( \frac{p_i}{5} )=1$,  $ ( \frac{d}{5}  )=1$, $X^d(5)(\Q_{5}) \neq \emptyset$.

	For all other primes $q$, $X^d(5)(\Q_{q}) \neq \emptyset$ by first and second parts of Theorem \ref{thm:ALL}.

	The case $N=13$ is quite similar to the case $N=5$, since corresponding Hilbert class field is $\Q(\sqrt{13},i)$ and they are both $1$ mod $4$.
	
	\item $N=6$: Hilbert class field of $\Q(\sqrt{-6})$ is $\Q(\sqrt{-3},\sqrt{2})$ and $\Q(j(\sqrt{-6}))$ is $\Q(\sqrt{2})$. Therefore using Theorem $\ref{thm:ALL}$ part 5, $X^d(6)$ has $\Q_{p_i}$-rational points if and only if $ ( \frac{2}{p_i} )=1$. 
	
	If $3$ is inert in $\K$ or splits in $\K$ then $X^d(6)(\Q_{3}) \neq \emptyset$ by parts 1 and 3 of Theorem \ref{thm:ALL}.

	If $2$ splits in $\K$ then $X^d(6)(\Q_{2}) \neq \emptyset$ by part 1 of Theorem \ref{thm:ALL}.

	And $2$ can not be inert in $\K$, since $ ( \frac{2}{p_i}  )=1$ for all $p_i$. 
	
	For all other primes $q$, $X^d(6)(\Q_{q}) \neq \emptyset$ by first and second parts of Theorem \ref{thm:ALL}.
	\item $N=10$: Hilbert class field of $\Q(\sqrt{-10})$ is $\Q(\sqrt{-2},\sqrt{5})$ and $\Q(j(\sqrt{-10}))$ is $\Q(\sqrt{5})$. 
	
	Therefore using Theorem $\ref{thm:ALL}$ part 5, $X^d(10)(\Q_{p_i}) \neq \emptyset$ if and only if $ ( \frac{5}{p_i} )=1$.

	Since $ ( \frac{5}{p_i} )= ( \frac{p_i}{5} )=1$,  $ ( \frac{d}{5} )=1$, $X^d(10)(\Q_{5}) \neq \emptyset$.

	If $2$ is inert or $2$ splits in $\K$ then $X^d(10)(\Q_{2}) \neq \emptyset$ by parts 1 and 4 of Theorem \ref{thm:ALL}.

	For all other primes $q$, $X^d(10)(\Q_{q}) \neq \emptyset$ by first and second parts of Theorem \ref{thm:ALL}.
\end{itemize}
\end{proof}

   Another result along these lines, which is a necessary condition for the existence of degree-$N$ $\Q$-curves, was given in \cite{Quer}:
 \begin{theorem} \label{q} \emph{(\cite{Quer}, Theorem 6.2)} Assume that there exists a quadratic $\Q$-curve of degree $N$ defined over some quadratic field $K$. Then every divisor $N_1 |N$ such that $$N_1 \equiv 1 (mod \;\;4)\; or\; N_1\;\; is \;\;even\;\; and \;\;N/N_1 \equiv 3 (mod\;\; 4)$$ is a norm of the field $K$.
 \end{theorem}
 
 Proof of this theorem is analytic, by constructing some functions on $X_0(N)$ with rational Fourier coefficients and studying the action of the involution $w_N$ on them. We'll take a more algebraic approach and given any square-free, relatively prime $d$ and $N$ show that Theorem \ref{thm:ALL} implies Theorem \ref{q} in the following two corollaries.

 Recall that saying that `$N_1$ is a norm in $\K$' is equivalent to say that $(N_1,d)=1$ where $(-,-)$ denotes the Hilbert symbol. Moreover $(N_1,d)=1$ if and only if all local Hilbert symbols, $(N_1,d)_p=1$ for all primes $p$. In other words, $(N_1,d)=-1$ if and only if there is a prime $p$ such that $(N_1,d)_p=-1$. Therefore, Theorem $\ref{q}$ gives a condition on the existence of local points. 
 
 The local Hilbert symbol is given by the following formula \cite{serreACA}:
  
  \begin{itemize}
	\item $p=\infty$: $(a,b)_{\infty}=1$ if and only if $a$ or $b$ is positive.
	\item $p\neq 2$: $(a,b)_{p}=(-1)^{\alpha\beta\epsilon(p)}\left( \frac{u}{p}\right)^\beta \left(\frac{v}{p}\right )^\alpha$ 
	\item $p=2$: $(a,b)_2=(-1)^{\epsilon(u)\epsilon(v)+\alpha w(v)+\beta w(u)}$ 
\end{itemize}
where $a=up^\alpha, b=vp^\beta$, $u,v$ are $p$-adic units in $\Q_p$ and $\epsilon(u), w(u)$ denote the class modulo $2$ of $\frac{u-1}{2}$ and of $\frac{u^2-1}{8}$, respectively. Therefore if $(N_1,d)_p=-1$ for some prime $p$ then $p$ divides $N_1$ or $d$. Since $\Pi_p (a,b)_p=1$, one can deduce that $(N_1,d)_p=-1$ for some odd prime $p$ that is dividing $N_1$ or $d$.
 
 
 \begin{corollary}\label{cor:odd}
 Let $N$ be an odd square-free integer such that there exists a divisor $N_1$ of $N$, $N_1 \equiv 1 (mod \;\;4)$ and $(N_1,d)_p=-1$ for some $p$. Then  $X^d(N)(\Q_p)=\emptyset$.
 \end{corollary}
 
  \begin{proof}
 Say $p | N_1$. Since $(N_1,d)_p=\left( \frac{d}{p} \right )=-1$, $p$ is inert in $\K$. Since $N_1 \equiv 1$ mod $4$, either $p \equiv 1$ mod $4$ or $p \equiv 3$ mod $4$ and there is another divisor $p'$ of $N_1$ that is also congruent to $3$ mod $4$. If $p \equiv 1$ mod $4$, then $X^d(N)(\Q_p)=\emptyset$ and if $p \equiv 3$ mod $4$, then there are at least $2$ primes dividing $N_1$ that are congruent to $3$ mod $4$, hence $X^d(N)(\Q_p)=\emptyset$, by part 3 of Theorem \ref{thm:ALL} .
 
 Say $p |d$. Since $(N_1,d)_p=\left( \frac{N_1}{p} \right )=-1$, $p$ is inert in $\Q(\sqrt{N_1})$. Let $H_O$ denote the ring class field of the order $\Z[\sqrt{-N}]$. Then $H_0=\Q(\sqrt{-N}, j(\sqrt{-N}))$ and $H_O \cap \mathbb{R}=\Q(j(\sqrt{-N}))$ by class field theory. Since $N_1 \equiv 1$ mod $4$, $\Q(\sqrt{N_1})$ lies in the genus field of $\Q(\sqrt{-N})$, hence $\Q(\sqrt{N_1}) \subset \Q(j(\sqrt{-N}))$. This shows that there is no prime of $\Q(j(\sqrt{-N}))$ lying above $p$ with residue degree $1$, thus $X^d(N)(\Q_p)=\emptyset$ by part 5 of Theorem \ref{thm:ALL}. 
  \end{proof}
 
 \begin{corollary}\label{cor:even}
 Let $N$ be an even square-free integer such that there exists an even divisor $N_1$ of $N$, $N/N_1 \equiv 1 (mod \;\;4)$ and $(N_1,d)_p=-1$ for some $p$. Then  $X^d(N)(\Q_p)=\emptyset$.
 \end{corollary}
   \begin{proof}
   Say $p|N_1$. Since $(N_1,d)_p=\left( \frac{d}{p} \right )=-1$, $p$ is inert in $\K$ and since $\prod_\nu (a,b)_\nu=1$, we can assume that $p$ is odd. By part 3 of Theorem \ref{thm:ALL}, $X^d(N)(\Q_p)\neq \emptyset$ if and only if $N=2p\Pi_i q_i$ such that $p$ is $3 $ mod $4$ and $q_i$ is $1$ mod $4$. Therefore $N_1=2p\Pi_i q_i$ for some $q_i$ congruent to $1$ mod $4$, this contradicts to the assumption that $N/N_1 \equiv 3$ mod $4$, hence $X^d(N)(\Q_p)=\emptyset$. 
   
   Say $p|d$, then $(N_1,d)_p=\left( \frac{N_1}{p} \right )=-1$. Again we can assume $p$ is odd. Let $H_O$ denote the ring class field of the order $\Z[\sqrt{-N}]$. Then $H_0=\Q(\sqrt{-N}, j(\sqrt{-N}))$ and $H_O \cap \mathbb{R}=\Q(j(\sqrt{-N}))$ by class field theory. If $N$ has odd number of prime factors that is congruent to $3$ mod $4$ then $N_1$ has even number of prime factors that is congruent to $3$ mod $4$ and genus field of $\Q(\sqrt{-N})$ includes $\sqrt{2},\sqrt{-l_i}, \sqrt{q_i}$ where $l_i$ and $q_i$ are prime factors of $N$ that are $3$ and $1$ mod $4$ respectively. If $N$ has even number of prime factors that is congruent to $3$ mod $4$ then $N_1$ has odd number of prime factors that is congruent to $3$ mod $4$ and genus field of $\Q(\sqrt{-N})$ includes $\sqrt{-2},\sqrt{-l_i}, \sqrt{q_i}$. In either case, the genus field includes $\Q(\sqrt{N_1})$. And since $H_O \cap \mathbb{R}=\Q(j(\sqrt{-N}))$, $\Q(\sqrt{N_1})$ lies in $\Q(j(\sqrt{-N}))$, where $H_O$ is ring class field of $\Z[\sqrt{-N}]$. Hence every prime of $\Q(\sqrt{-N})$ that is lying above $p$ has a residue degree greater than $1$, therefore $X^d(N)(\Q_p)=\emptyset$ by part 5 of Theorem \ref{thm:ALL}.   
   \end{proof}

Corollaries \ref{cor:odd} and \ref{cor:even} imply Theorem \ref{q}.

Another result about existence of local points on $X^d(N)$ is given by Clark in \cite{Clark}. Generalizing the techniques that is used in this proof we prove part 3 of Theorem \ref{thm:ALL}. As a result, the following theorem follows from part 3 of Theorem \ref{thm:ALL}. 

\begin{theorem}\emph{(\cite{Clark})}\label{thm:Clark}Let $N$ be a prime number congruent to $1$ mod $4$, $d=p^*=(-1)^{\frac{p-1}{2}}p$ where $p$ is a prime different than $N$ and $\left ( \frac{N}{p}  \right )=-1$. Then $X^d(N)(\Q_N)=\emptyset$.
\end{theorem}

In the same paper (\cite{Clark}), it was asked whether or not $p$ and $N$ were the only primes that $X^d(N)$ fails to have local points. We prove that the answer is `yes' in Corollary \ref{cor:Clark}.

\textbf{Corollary:}(Corollary \ref{cor:Clark}) Let $N$ be a prime congruent to $1$ mod $4$, $p$ be an odd prime such that $(N/p)=-1$ and $d=p^*$ then 
\begin{enumerate}
	\item $X^d(N)(\Q_N)=\emptyset$
	\item $X^d(N)(\Q_p)=\emptyset$
	\item $X^d(N)(\Q_\ell) \neq \emptyset$ for any other prime $\ell$ different than $p$ and $N$.
\end{enumerate}
 
   \section{Real Points}\label{sec:real}
 
 We will keep the same notation as in the previous sections. Given a square-free integer $N$, a quadratic number field $\K:=\Q(\sqrt{d})$ and a prime $p$ what can be said about $X^d(N)(\Q_p)$, where $X^d(N)$ is the twist of $X_0(N)$ with $w_N$ and $<\sigma>:=\Gal(\Q(\sqrt{d})/\Q)$. Since we are dealing with local points, by abuse of notation we regard $\sigma$ as the generator of the extension $\K_{\nu}$ over $\Q_p$, where $\nu$ is a prime of $\K$ lying over $p$. Let $k$ be the residue field and $R$ be the valuation ring of $\K_{\nu}$.
 
 We start with the real points of $X^d(N)$:
 
 \begin{proposition} $X^d(N)(\R) \neq \emptyset$. \label{propreal}
\end{proposition}

\begin{proof}
If $\K$ is an imaginary quadratic field, $\K_{\nu}=\C$, $\Q_{\infty}=\R$ and $\sigma$ is complex conjugation.
Say $E$ has CM with the full ring of integers of $\Q(\sqrt{-N})$. Then $E$ corresponds to the
lattice (a) $ [1, \sqrt{-N}]$ or (b) $[1,(1+\sqrt{-N})/2]$ and $E$ induces a fixed point, $P$ of $w_N$, in
$X_0(N)(\C)$.

Now for (a), conjugate of  $[1, \sqrt{-N}]$  is $[1, -\sqrt{-N}]$. And the matrix
$\left ( \begin{array}{cc}
	1 & 0 \\
	0 & -1 \\
\end{array} \right )$
sends the first basis to the second.

For (b), conjugate of $[1,(1+\sqrt{-N})/2]$ is $[1,(1-\sqrt{-N})/2]$ and the matrix
$\left ( \begin{array}{cc}
	1 & 0 \\
	1 & -1\\
\end{array} \right )$
sends the first basis to the second.

 Since these matrices are invertible
in $\Z$, they both induce the same elliptic curve i.e. $P$ is a real point of $X^d(N)$ that is fixed by $w_N$, hence induces a point of $X^d(N)(\R)$.

If $\K$ is real quadratic then $\sigma$ induces trivial map and same argument above shows that $X^d(N)(\R) \neq \emptyset$.
\end{proof}

 Now lets assume $p$ splits in $\K$ then a copy of $\K$ is in $\Q_p$. Since $X_0(N)$ and $X^d(N)$ are isomorphic over $\K$ and $X_0(N)(\Q_p)$ is non-empty, $X^d(N)(\Q_p)$ is non-empty.
 
 Therefore, $X^d(N)$ might fail to have $p$-adic points only for finite primes that are inert or ramified in $\K$.
 
 \section{Primes That Are Inert in $\Q(\sqrt{-N})$}\label{sec:inert}
 




For the two main cases that we are dealing-the inert case and ramified case- we will be using different tools. For the inert case, since we have the notion of etale descent for $X_0(N)_{/ R}$, existence or non-existence of local points will be showed by checking the existence of points over the corresponding special fiber. For this, the following version of Hensel's Lemma  will be used:

\begin{lemma}\label{lem:hensel}Let $K$ be a complete local ring , $R$ its valuation ring and $k$ its residue field. Let $X$ be a regular scheme over $S:=Spec(R)$ and $f:X \rightarrow S$ a proper flat morphism. Say $X_{\eta}:=X \times_S Spec(K)$ is the generic fiber and $X_o:=X \times_S Spec(k)$, the special fiber. Then each $K$-rational point of $X_{\eta}$ corresponds to a \emph{smooth} $k$-rational point of $X_o$. 
\end{lemma}

 We will also use the following theorems of Deuring:

\begin{theorem} \label{thm:Deuring} Let 
\begin{itemize}
	\item $\tilde{E}$ be an elliptic curve over a number field $L$ that has CM by $\Q(\sqrt{-N})$;
	\item $p$ be a rational prime;
	\item $\beta$ be a prime of $L$ lying over $p$ such that $\tilde{E}$ has good reduction over $\beta$ and
	\item $E$ denote the reduction of $\tilde{E}$.
\end{itemize}
    Then $E$ is supersingular if and only if $p$ is ramified or inert in $\Q(\sqrt{-N})$.
\end{theorem} 
   
 \begin{theorem}\label{thm:DeuringLifting}\emph{(Deuring Lifting Theorem, see \cite{Lang} Theorem 14) }
 Let $E$ be an elliptic curve over a finite field $k$ of characteristic $p$, and let $\alpha$ be an element of $End(E)$. Then there exists:
 \begin{itemize}
	\item an elliptic curve $\tilde{E}$ over a number field $B$
	\item an endomorphism $\tilde{\alpha} \in End(\tilde{E})$ and
	\item a place $\beta$ of $B$ lying over $p$
\end{itemize}
such that 
\begin{itemize}
	\item the reduction of $\tilde{E}$ mod $\beta$ is $E$
	\item the reduction of $\tilde{\alpha}$ mod $\beta$ is $\alpha$ and
	\item $|k|=p^f$ where $f$ is the inertia degree of $\beta$ over $p$.
\end{itemize}
\end{theorem} 

When $p$ is inert in $\Q(\sqrt{d})$, $\sigma$ induces the non-trivial map(Frobenius) on $\Gal(k/\F_p)$, where $k=\F_{p^2}$. We have different cases according to $p |N$ or not. 

\subsection{$p$ dividing the level}

Say $p \mid N$, $\nu$ the prime of $\K$ lying over $p$ and $R$ is the ring of integers of the localization $\K_{\nu}$.

\begin{figure}[h]
	\centering
		\includegraphics{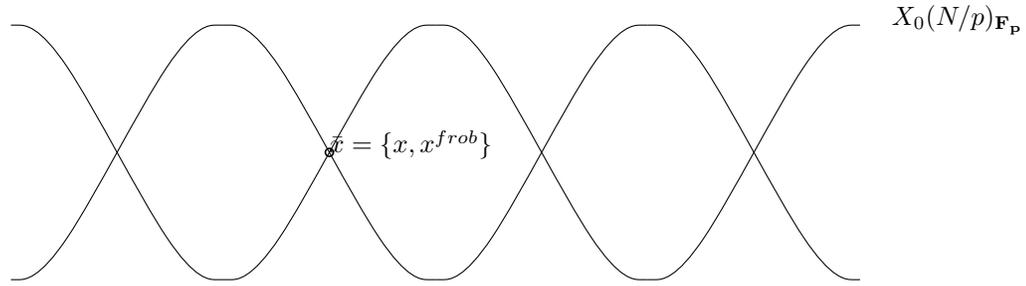}
	\caption{Special fiber of $\mathcal{X}_0(N)/_{\Z_p}$}
	\label{fig:model2}
\end{figure}

 In Mazur\cite{Mazur2} and Deligne-Rapaport\cite{Deligne} there is a model of $\mathcal{X}_0(N)_{/\Z_p}$ whose special fiber $X_0(N)_{\F_p}$ is 2 copies of $X_0(N/p)_{\F_p}$ glued along supersingular points twisted by first power Frobenius. The Atkin-Lehner involution $w_N$ interchanges the two branches and Frobenius stabilizes each branch, see Figure \ref{fig:model2}. A regular model $\tilde{\mathcal{X}}_0(N)$ can be obtained by blowing-up ($|Aut(E,C)|/2-1$)-many times at each supersingular point. Actions of Frobenius and $w_N$ extends to regularization as well.

Since $p$ is ramified in $\Q(\sqrt{-N})$, by Theorem \ref{thm:Deuring}, any $w_N$-fixed point is supersingular. Say $x$ is a $w_N$-fixed supersingular(hence singular) point of $X_0(N)_{/ \F_p}$. To have a regular model, we need to blow-up at each supersingular point ($|Aut(E,C)|/2-1$)- times where $(E,C)$ is on $X_0(N/p)_{/\F_p}$. To keep track of the different schemes, we need to introduce some notation. As introduced in the beginning, $\mathcal{X}_0(N)$ denotes the model of $X_0(N)$ over $\Z_p$(not necessarily regular). Let $\tilde{\mathcal{X}}_0(N)$ denote the regularization of $\mathcal{X}_0(N)$ after blow-ups, and $\tilde{X}_0(N)_{/\F_p}$ be its special fiber.

We can define a model of $X^d(N)$ over $\Z_p$ as a descent to $Spec \Z_p$ of $X_0(N) x_{Spec \Z_p} Spec R$ by a descent datum twisted by $w_N$
Note that extension $R/\Z_p$ is Galois since $p$ is inert in $\K$. This model will be denoted as $\mathcal{X}^d(N)_{/\Z_p}$. Our aim is to use Hensel's Lemma(Lemma \ref{lem:hensel}) to make conclusions about $\Q_p$-rational points of the generic fiber of $\mathcal{X}^d(N)_{/Z_p}$. In order to do this, we must first show that $\mathcal{X}^d(N)_{/\Z_p}$ is regular.

\begin{proposition}\label{prop:flat}
Let $R\hookrightarrow S$ be a flat extension of local rings. If $S$ is regular than so is $R$.
\end{proposition}

\begin{proof}
Matsumura, Commutative Algebra, page 155, Theorem 51.
\end{proof}

\begin{proposition}
$\mathcal{X}^d(N)_{/\Z_p}$ is a regular model of $X^d(N)$.
\end{proposition}

\begin{proof}
By definition of etale descent, $\mathcal{X}^d(N)$ and $\tilde{\mathcal{X}}_0(N)$ are isomorphic over $R$. Since $R$ is an unramified extension of $\Z_p$, $\tilde{\mathcal{X}}_0(N) \times _{\Z_p} R$ is regular, hence $\mathcal{X}^d(N) \times _{\Z_p} R$ is also regular. By Proposition \ref{prop:flat}, $\mathcal{X}^d(N)_{/\Z_p}$ is regular.
\end{proof}


Now we will give a necessary condition for the existence of a smooth point on $X_0^d(N)(\F_p)$.

\begin{proposition}\label{prop:degree}
 There exists a smooth point on $X^d(N)(\F_p)$ if and only if there is a point on $X_0(N/p)(\F_{p^2})$ corresponding to a supersingular elliptic curve with a degree-$4$ automorphism.
\end{proposition}

\begin{figure}[h]
	\centering
		\includegraphics{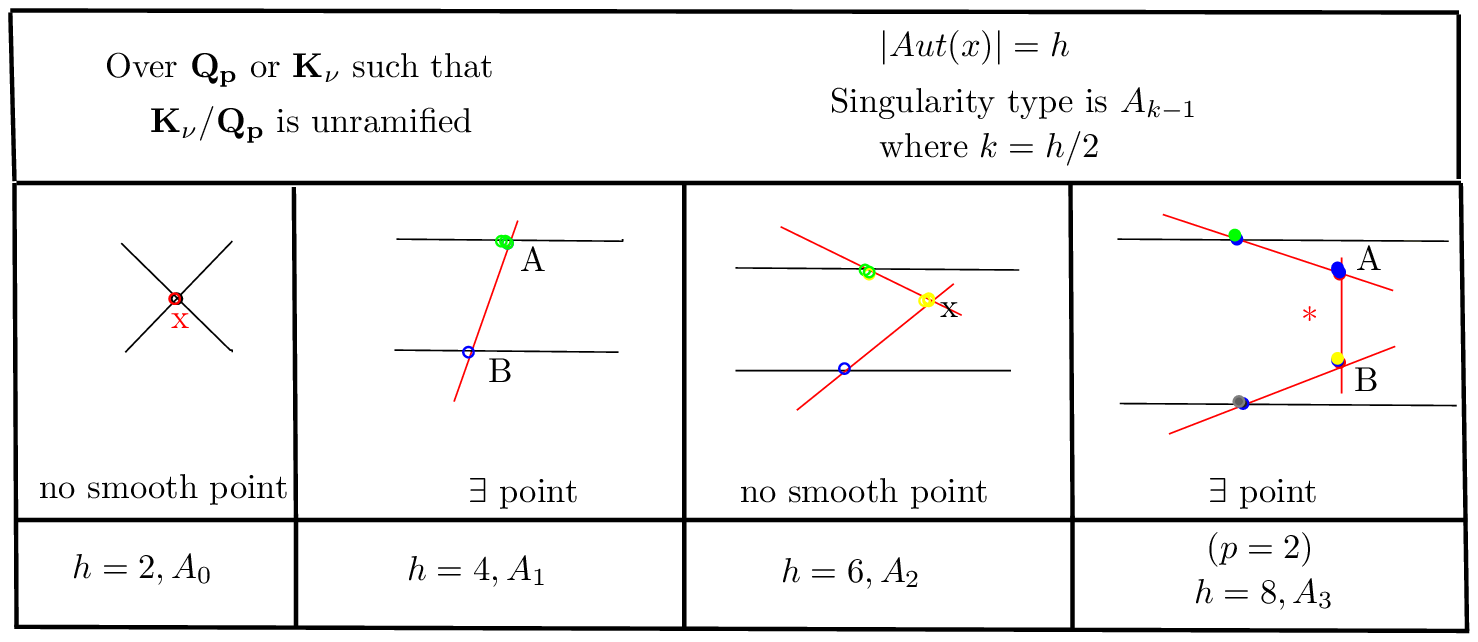}
	\caption{Blow-ups}
	\label{fig:TABLE}
\end{figure}

\begin{proof}
As explained above, $X^d(N)$ is generic fiber of $\mathcal{X}^d(N)_{/\Z_p}$ which is etale descent of $\mathcal{X}_0(N)_{/\Z_p}$ by $R/\Z_p$. Since $w_N$ interchanges the branches of $X_0(N)_{/\F_p}$, $\F_p$-rational points on the special fiber $X^d(N)_{/\F_p}$ are coming from supersingular points of $X_0(N)_{/\F_p}$, which are all singular. And we must keep in mind that at each singular(hence supersingular) point we have $(|Aut(E,C)|/2-1)$-many exceptional lines. Note that automorphism group of an elliptic curve over a field of characteristic $\ell$ is $\mu_2,\mu_4$ or $\mu_6$ if $\ell$ is not $2$ or $3$ where $\mu_s$  denotes the group of primitive $s$-th roots of unity. If $\ell=2$ or $3$ and $E$ is the unique supersingular elliptic curve in characteristic $\ell$ then $Aut(E)$ is $C_3 \rtimes \{\pm1,\pm i,\pm j,\pm k\}$ or $C_3 \rtimes C_4 $ respectively, where $C_m$ denotes the cyclic group of order $m$. 


Therefore if $|Aut(P)|=4n$ for $n>1$, there is an element of order $4$ in $Aut(P)$ and the number of blow-ups is $2n-1$ which is odd. Since we have odd number of exceptional lines, then there is one line $L_{/ \F_p}$ that is fixed by the action of $w_N$(see the second column of Figure \ref{fig:TABLE}). On this line $L$ the points $A$ and $B$ are singular and fixed by $w_N \circ \sigma$ but these are not the only fixed points. The action of $\sigma \circ w_N$ on zeroth, first and second cohomology of $L$ has traces $1,0$ and $p$ respectively. Then by Lefschetz fixed point theorem, there is a smooth $w_N \circ \sigma$-point on this exceptional line $L$. Therefore if we have a supersingular point with a degree-$4$ automorphism then there is a smooth point on $X^d(N)(\F_p)$.

For the reverse direction, say there is no such supersingular point $P$ with a degree-$4$ automorphism. If  $|Aut(P)|$ is $2$, then $\mathcal{X}_0(N)$ is already regular but all $\F_p$-rational points of $X^d(N)$ are singular. 

If $|Aut(P)|=6$, then we replace this point by 2 exceptional lines over $\F_p$ and $\sigma \circ w_N$ interchanges these lines. Each of these exceptional lines cuts one of the branches and also the other exceptional line once. Call the intersection point of these lines as $x$, then $\sigma(x)=x$ and it is the only point fixed by the action of $\sigma$ on these lines. Then $w_N(x)=w_N(\sigma(x))=\sigma(w_N(x))$, hence $w_N(x)$ is also fixed by $\sigma$ i.e. $w_N(x)=x=\sigma(x)$. Thus $x$ induces an $\F_p$-rational point of $X^d(N)$. However $x$ is a singular point. For a picture of this situation we refer to the table at the end of this section. 
\end{proof}

Using Proposition \ref{prop:degree} and Hensel's Lemma we get the following: 

\begin{corollary}\label{cor:degree}
There exists a point on $X^d(N)(\Q_p)$ if and only if there is a supersingular point with a degree-$4$ automorphism.
\end{corollary}

\begin{theorem}\label{thm:IDB}
Let $N$ be a square-free positive integer and $p$ be an odd prime that is inert in $\K$ and dividing $N$ then $X^d(N)(\Q_p) \neq \emptyset$ if and only if $p$ is $3$ mod $4$ and $N$ is of the form:
\begin{enumerate}
	\item $p \Pi_i q_i$ such that all $q_i \equiv 1$ mod $4$ or
	\item $2 p \Pi_i q_i$ such that all $q_i \equiv 1$ mod $4$.
\end{enumerate}
  
\end{theorem}

\begin{proof}By Corollary \ref{cor:degree} above, a local point exists if and only if there is a supersingular point $(E,C)$ on $X_0(N)_{/F_p}$ with automorphism group divisible by $4$. In particular we want $j=1728$ to be a supersingular $j$-invariant over characteristic $p$.

In order to have such a point, $\Z/4\Z$ must inject into $(\Z/q_i\Z)^*$ for every odd prime divisor $q_i$ of $N$ that is not $p$. Therefore if $N$ is odd, $N= p \Pi_i q_i$ such that $p \equiv 3$ mod $4$ and all $q_i \equiv 1$ mod $4$. 

Since the degree-$4$ automorphism $[i]$ sends a $2$-torsion point $(x,0)$ of $E_{1728}:y^2=x^3+x$ to $(-x,0)$, $[i]$ fixes the cyclic-$2$ subgroup $<(0,0)>$ of $E_{1728}[2]$. Hence if $N$ is even and there is a supersingular point on $X_0(N)_{/\F_p}$ with automorphism group divisible by $4$, then $N=2 p \Pi_i q_i$ such that $p \equiv 3$ mod $4$ and all $q_i \equiv 1$ mod $4$.

Conversely, if $N$ is one of the given form then $1728$ is a supersingular $j$-invariant. Let $E_{1728}$ be the elliptic curve over $\F_p$ of which $j$ invariant is $1728$. Then $[i]$ is in $Aut(E_{1728})$ and acts on $E_{1728}[\Pi_i q_i]$. The automorphism $[i]$ stabilizes a cyclic-$\Pi_i q_i$ subgroup if and only if $[i]$ stabilizes cyclic-$q_i$ subgroups of $E_{1728}[q_i]=\Z/q_i\Z \times \Z/q_i\Z$ for all $i$. The automorphism $[i]$ can be seen as an element of $GL_2(\F_{q_i})$ and it stabilizes a cyclic subgroup of order $q_i$ if and only if $[i]$ has eigenvalues defined over $\F_p$. If $q_i$ is odd then minimal polynomial of $[i]$ is $x^2+1$, this is equivalent to say that $q_i \equiv 1(mod\;\; 4)$ for all $i$. If $q_i=2$ then minimal polynomial of $[i]$ is $x+1$ and $[i]$ fixes the cyclic-$2$ subgroup $<(0,0)>$ of $E_{1728}[2]$.  
\end{proof}

For $p=2$ inert in $\K$ and $N$ even we get the following result: 

\begin{theorem}\label{thm:inerteven} If $2$ is inert in $\K$ and dividing $N$ then $X^d(N)(\Q_2)\neq \emptyset$ if and only if $N=2\Pi_i q_i$ such that $q_i \equiv 1$ mod $4$.
\end{theorem}

\begin{proof}
Over $\F_2$, $1728$ is the only supersingular $j$-invariant and $|Aut(E_{1728})|=24$. Say $q_i$ is a prime dividing $N/2$. Since $N$ is square-free, $q_i$ is odd. 

By Corollary \ref{cor:degree} above, a $\Q_2$-point exists if and only if there is a supersingular point $(E_{1728},C)$ on $X_0(N)_{/\F_2}$ with automorphism group divisible by $4$. In order to have such a point, $\Z/4Z$ must inject into $(\Z/q_i\Z)^*$ for every odd prime divisor $q_i$ of $N/2$ , hence $N= 2 \Pi_i q_i$ such that $q_i \equiv 1$ mod $4$ for all $i$.

For the converse, the argument is the same as in the corresponding part of Theorem \ref{thm:IDB}.

\end{proof}


\subsection{$p$ does not divide $N$}
We will be using the same notations introduced in the previous subsection, in particular $\mathcal{X}^d(N)_{/\Z_p}$ denotes etale descent $\mathcal{X}_0(N)$ from $R$ to $\Z_p$. If $p$ is not dividing $N$ then the following models are smooth, in particular regular:
\begin{itemize}
	\item $\mathcal{X}_0(N)_{/\Z_p}$,
	\item $\mathcal{X}_0(N) \times_{\Z_p} R$ (since $R/\Z_p$ is unramified) and
	\item $\mathcal{X}'^d(N)_{/\R}:=\mathcal{X}^d(N) \times_{\Z_p} R$ (since $\mathcal{X}'^d(N)_{/\R}$ is isomorphic to $\mathcal{X}_0(N) \times_{\Z_p} R$)
\end{itemize}

  By Proposition \ref{prop:flat} $\mathcal{X}^d(N)_{/\Z_p}$ is also regular. We will construct a point on the special fiber $X^d(N)(\F_p)$ using theory of CM curves and then by Hensel's Lemma we will be done.

  Let $\Sigma_N$ be the set of tuples $(E,C)$ such that $E$ is a supersingular elliptic curve over characteristic $p$ and $C$ is cyclic group of order $N$. We start by studying the action of the involution $w_N \circ \sigma$ on $\Sigma_N$.

\begin{definition} Let $B$ be the unique quaternion algebra over $\Q$ that is ramified only at $p$ and infinity. An order $O$ of $B$ is of level $N$ if 
\begin{itemize}
\item $q \neq p$,  $O_q = O \otimes_{\Z}\Z_q \cong 
\left ( \begin{array}{cc}
	\Z_q & \Z_q \\
	N\Z_q & \Z_q \\
\end{array} \right )$
\item $O_p \cong \{ \left ( \begin{array}{cc}
	\alpha & \beta \\
	p\bar{\beta} & \bar{\alpha} \\
\end{array} \right )| \alpha, \beta \in R \}$ where $R$ is the ring of integers of the unique unramified quadratic extension of $\Q_p$.
\end{itemize}
\end{definition}

\begin{proposition}
Let $B:=End(E) \otimes_{\Z} \Q$, $E$ is a supersingular elliptic curve over $\F_{p^2}$ and $C$ a cyclic subgroup of $E$ of order $N$. Then $B$ is the unique quaternion algebra over $\Q$ which is ramified only at $p$ and $\infty$, $End(E)$ is a maximal order in $B$ and $End(E,C)$ is an Eichler order of level $N$.
\end{proposition}

\begin{proof}
For $B$ being the claimed quaternion algebra and $End(E)$ its maximal order we refer to Silverman \cite{Silverman1} Chapter 5. Assuming this let $(E,C)$ be an element of the set $\Sigma_N$. The map $\pi: E \rightarrow E/C$ induces an isomorphism between $B = End(E) \otimes \Q$ and $End(E/C) \otimes \Q$ via the map $f \mapsto \pi f \pi^{-1}$. Since $End(E/C)$ is a maximal order in the latter, its image is a maximal order in the former. Moreover the intersection of these two maximal orders is $End(E,C)$, hence an Eichler order.  
\end{proof}

It's important to understand what kind of maps are in $R=End(E,C)$. Since we have an algebraic way to study these rings, properties of endomorphisms of $(E,C)$ can be also understood  by studying orders embedded in $R$.

\begin{definition}

Let $L$ be an imaginary quadratic number field, $O$ an order of $L$, $\alpha:L \hookrightarrow B$ an algebra embedding such that $\alpha(L) \cap R =\alpha(O)$ where $R$ is an Eichler order of level $N$ in $B$ as above. Then the pair $(R,\alpha)$ is called as an \emph{optimal embedding of} $O$.

\end{definition}

The following theorem of Eichler states conditions for existence of an optimal embedding:

\begin{proposition}\label{Prop:Eic}

Given an $R$, $B$ as above and $L=\Q(\sqrt{M})$, an optimal embedding $(R,\alpha)$ of an order $O$ of $L$, exists if and only if 
\begin{itemize}
\item $M <0$, $p$ is inert or ramified in $L$ and $p$ is relatively prime to the conductor of $O$ and
\item $q$ splits or ramifies in $O$ for every $q$ dividing $N$.
\end{itemize}

\end{proposition}

For any $q | N$ and $q':=N/q$, let $w_q$ be the Atkin-Lehner operator that sends $(E,C) \mapsto (E/q'C, E[q]+C/q'C)$ where $E[q]$ is the kernel of multiplication by $q$. Each $w_i$ is an involution and $w_i \circ w_j=w_j \circ w_i=w_{ij}$ for every coprime $i,j$. 
We have another operator, \emph{Frobenius} acting on the set $\Sigma_N$ and remember that $\sigma$ is acting as Frobenius in the inert case. The following result shows that Frobenius also can be seen as an Atkin-Lehner operator on the set $\Sigma_N$.

\begin{theorem}(see Chapter V, Section 1 of \cite{Deligne})\label{thm:ribetfrob}
The involution $w_p$ permutes the two components of $X_0(Np)_{\F_p}$. It acts on the set of singular points of $X_0(Np)_{\F_p}$ as the Frobenius morphism $x \mapsto x^p$. 
\end{theorem}

Let $\psi$ be a map from $X_0(N)_{/ \F_p}$ to $X_0(Np)_{/ \F_p}$, an isomorphism onto one of the two components. The map $\psi$ takes the supersingular locus of $X_0(N)_{/ \F_p}$ to the supersingular locus of $X_0(Np)_{/ \F_p}$. The set $\Sigma_N$ defined in the beginning of the section is the supersingular locus of $X_0(N)_{/ \F_p}$. The Atkin-Lehner operator $w_{Np}$ acts on $X_0(Np)$, in particular acts on $\psi(\Sigma_N)$. When we say the action of $w_{Np}$ on $\Sigma_N$, we actually mean the action of $w_{Np}$ on $\psi(\Sigma_N)$.

\begin{corollary}\label{cor:optimal}
Given $B$ there is an embedding of $\Z[\sqrt{-pN}]$ into some Eichler order $R$ of level $N$.
\end{corollary}

\begin{proof}
Using the first part of Proposition \ref{Prop:Eic} quadratic imaginary field $L=\Q(\sqrt{-pN})$ embeds in $B$ since $p$ is ramified in $\Q(\sqrt{-pN})$ and the embedding is optimal for the maximal order of $\Q(\sqrt{-pN})$ since any prime $q$ dividing $N$ is ramified in $\Z[\sqrt{-pN}]$. If the maximal order $O_L$ of $L$ is $\Z[\sqrt{-pN}]$ we are done. Say $\Z[\sqrt{-pN}]$ is of conductor $2$ in $O_L$ then $\alpha(\Z[\sqrt{-pN}])\subset \alpha(O_L)=\alpha(L) \cap R  $, hence $\Z[\sqrt{-pN}]$ embeds in $R$.
\end{proof}

\begin{corollary}\label{cor:wNp}
If there is an embedding $\alpha$ of $\Z[\sqrt{-pN}]$ into $R$ for some Eichler order $R=End(E,C)$ of level $N$ then there is a fixed point of $w_{Np}$ in $\Sigma_N$.
\end{corollary}

\begin{proof}
By assumption there exists an element whose square is $-pN$ in $R$ i.e. an endomorphism of degree $Np$ of $(E,C)$, in particular an endomorphism of degree $Np$, say $f$, of $E$. Using Deuring's Lifting Theorem(Theorem \ref{thm:DeuringLifting}) this endomorphism and $E$ can be lifted to characteristic $0$ i.e. $(\tilde{E},ker(\tilde{f}))$ ($E$ and $f$ lifted to characteristic $0$) is in $X_0(Np)(\bar{\Q})$ and is fixed by $w_{Np}$. The reduction of this point mod $p$,is a supersingular point on $X_0(Np)_{/ \F_p}$ that is fixed by $w_{Np}$ since $E$ has no $p$-torsion $|ker(f)|=N$, $(E,ker(f))$ is identified with a point in $\Sigma_N$ i.e. is in $\psi(\Sigma_N)$. 
\end{proof}

\begin{example}\label{ex:noncm}
Let $p=7$ and $N=5$. Since $\left ( \frac{-5}{7}  \right )=1$,by Theorem \ref{thm:Deuring} reduction of any elliptic curve which has CM by $\Q(\sqrt{-5})$ over $p$ is ordinary. Hence there isn't any $w_5$-fixed point in $\Sigma_5$ i.e. there is no optimal embedding of $\Z[\sqrt{-5}]$ into any $R$ where $R$ is an Eichler order of level $5$ in the quaternion algebra $\Q_{7,\infty}$. In fact, there isn't any embedding of $\Q(\sqrt{-5})$ into $\Q_{7,\infty}$ since $7$ splits in $\Q(\sqrt{-5})$ the localization of $\Q(\sqrt{-5})$ at the primes lying above $7$, is not even a field.
\end{example}

\begin{theorem}\label{thm:inertnotdividing}
If $p$ is inert in $\Q(\sqrt{d})$ and $p \nmid N$ then $X^d(N)(\Q_p) \neq \emptyset$.
\end{theorem}

\begin{proof}
By Corollaries \ref{cor:optimal} and \ref{cor:wNp},  there is a point $x \in \Sigma_N$ such that $w_{Np}(x)=x$. Since $(p,N)=1$ we have $w_N \circ w_p(x)=w_{Np}(x)=x$. By Theorem \ref{thm:ribetfrob}, $w_p$ acts as $frob_p$ on the set $\Sigma_N$, hence $w_N \circ frob_p(x)=x$. By theory of etale descent this gives a point in $X^d(N)(\F_p)$ and since $p \nmid N$, we have smooth model, by Hensel's Lemma(Lemma \ref{lem:hensel}) we are done.
\end{proof}

\section{Primes Ramified in $\K$ and Unramified in $\Q(\sqrt{-N})$}\label{sec:ram}

Let $p$ be a prime that is ramified in the quadratic field $\K$ but not in $\Q(\sqrt{-N})$. Let $\nu$ be the prime of $\K$ lying over $p$, and $R$ be the ring of integers of $\K_{\nu}$. Note that we don't have a good model for $X^d(N)$ over $\R$, since $\R/\Z_p$ is not etale.
By assumption, $p \nmid N$. Then by a well-known theorem of Igusa, $\mathcal{X}_0(N)$ is a smooth $\Z[1/N]-$scheme, hence for any $p \nmid N$ the special fiber of $\mathcal{X}_0(N) \rightarrow Spec(R)$ is \emph{smooth} over the residue field $\R/\nu$.
 
 Since $p$ is ramified, the residue field $\R/\nu$ is $\F_p$, and the induced action of $\sigma$ on the residue field is trivial. 
In  this setting, our approach will be to produce points on $X_0(N)(\Q_p)$ which are fixed by $w_N$ and which are thus CM points. Note that such points are $\Q_p$-rational points of $X^d(N)$. The main tool in showing the existence of such points is Deuring's Lifting Theorem which is stated as Theorem \ref{thm:DeuringLifting}. This theorem allows us to lift $w_N$-fixed points of $X_0(N)(\F_p)$ to $w_N$-fixed points of $X_0(N)(\Q_p)$, as Proposition \ref{prop:onlynice} demonstrates. Before stating the proposition we need to recall the following facts about CM elliptic curves.
 
  If $E$ corresponds to a fixed point of $w_N$ on $X_0(N)(\bar{\Q})$ then $E$ has an endomorphism whose square is $[-N]$. Hence $End(E)$ contains a copy of $\Z[\sqrt{-N}]$ and can be embedded in $\Z[\frac{D+\sqrt{D}}{2}]$ where $D$ is the discriminant of the CM field $\M:=\Q(\sqrt{-N})$. If $N \equiv 1$ or $2$ mod $4$ then these two orders are the same, hence $End(E)$ is the maximal order of $\M$. Otherwise, $End(E)$ is an order of conductor $2$ in the maximal order. Throughout the section we will use the following notation:
  
\begin{itemize}
  \item $O:=\Z[\sqrt{-N}]$ 
	\item $h$ denotes the class number of $O$
	\item $E$ is an elliptic curve such that $End(E)$ contains $\Z[\sqrt{-N}]$
	\item $\H_O$ is the ring class field of $O$ 
\end{itemize}
 
 Recall that by theory of CM, we have $h$-many elliptic curves which has CM by $O$ and their $j$-invariants are all conjugate. 
  
  \begin{proposition} \label{CM}
  Let $E$ be an elliptic curve over a number field $B$ and $E$ has an endomorphism $\alpha_0$ whose square is $[-N]$. Then $(E,ker(\alpha_0))$ is a $w_N$-fixed point on $X_0(N)(B)$.
  
  \end{proposition}
  
  \begin{proof}
  By definition $(E,ker(\alpha_0))$ is a $w_N$-fixed point of $X_0(N)(\bar{\Q})$ and $E$ is defined over $B$, $\alpha_0$ is defined over $B(\sqrt{-N})$.  Let $\phi$ be the generator of $Gal(B(\sqrt{-N})/B)$ then $ker(\alpha_0)^{\phi}=ker(\pm \alpha_0)=ker(\alpha_0)$ since the only endomorphisms of $E$ whose square is $[-N]$ are $\pm \alpha_0$. Therefore $ker(\alpha_0)$ is defined over $B$ as well. 
  \end{proof}

 \begin{proposition} \label{prop:onlynice}
 Let $p$ be an odd prime. Any $w_N$-fixed point on $X_0(N)(\F_p)$ is the reduction of a $w_N$-fixed $\Q_p$-rational point on the generic fiber of $X_0(N)$. Conversely, a $w_N$-fixed point on the generic fiber reduces to a $w_N$-fixed point on $X_0(\F_p)$.
 \end{proposition}
 
 \begin{proof}
 Say there is a $w_N$-fixed point on $X_0(N)(\F_p)$. This point corresponds to an elliptic curve $E$ over $\F_p$ which has an endomorphism $\alpha$ whose square is $[-N]$. By Theorem \ref{thm:DeuringLifting} we have the following:
 
\begin{itemize}
	\item An elliptic curve $\tilde{E}$ over a number field $B$,
	\item an endomorphism $\alpha_0$ whose square is $[-N]$ and 
	\item a prime $\beta$ over $p$ with residue degree $1$.
\end{itemize}
 
 such that 
 
\begin{itemize}
	\item Reduction of $\tilde{E}$ mod $\beta$ is $E$ and
	\item reduction of $\alpha_0$ mod $\beta$ is $\alpha$.
\end{itemize}
 
 By Proposition \ref{CM}, it is enough to check whether or not $B$ can be embedded in $\Q_p$. Since we have $f(\beta|p) = 1$ it suffices to check that $B$ is unramified at $p$. 
 
\begin{itemize}
	\item Say $O=\Z[\sqrt{-N}]$ is the maximal order of $\Q(\sqrt{-N})$ then $\H_O$ is the Hilbert class field hence $\H_O/\Q(\sqrt{-N})$ is unramified. By assumption $p$ is unramified in $\Q(\sqrt{-N})$ hence $p$ is unramified in $\H_O/\Q$. Therefore $p$ is unramified in $B/\Q$ since $B=\Q(j(E))$ lies inside $\H_O$.
	\item Say $O=\Z[\sqrt{-N}]$ is order of conductor $2$ in $O'=\Z[\frac{1+\sqrt{-N}}{2}]$. Any odd prime(that is not lying over $2$) is unramified in $\H_O/\Q(\sqrt{-N})$. Then the above argument shows that any $p>2$ is unramified in $B/\Q$. 
\end{itemize}
Therefore $e(\beta|p)=1$ and $B_{\beta} \cong \Q_p$. 
 Conversely, since the fixed locus of $w_N$ is proper, a $w_N$-fixed point on $X_0(N)(\Q_p)$ reduces to a $w_N$-fixed point on $X_0(N)(\F_p)$.
 \end{proof}

Proposition \ref{prop:onlynice} shows that if $X_0(N)(\F_p)$ contains a smooth $w_N$-fixed point, then $X^d(N)(\Q_p)$ is nonempty. We now show the converse.
 
 \begin{proposition}\label{prop:diag}
   Let $x$ be a point of $X_0(N)(\K_{\nu})$ such that $w_N(x^{\sigma})=x$ then $x$ reduces to a  $w_N$-fixed point on the special fiber of $\mathcal{X}_0(N)_{/R}$.
 \end{proposition}
 
 \begin{proof}
 Note that $\sigma$ is not a morphism of $Spec(R)$-schemes. We define the map $\hat{\sigma}:\mathcal{X}_0(N) \rightarrow \mathcal{X}_0(N)$ using the following diagram:
 
 
\begin{figure}[h]
	\centering
		\includegraphics{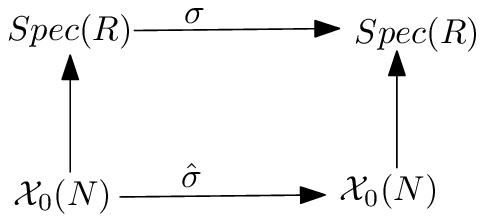}
	\label{fig:kucukdiagram}
\end{figure}

Since $\K_{\nu}/\Q_p$ is ramified, $\sigma$ induces the trivial action on the residue field $R/\nu$, therefore $\sigma$ induces trivial action on the special fiber: 

\begin{figure}[h]
	\centering
		\includegraphics{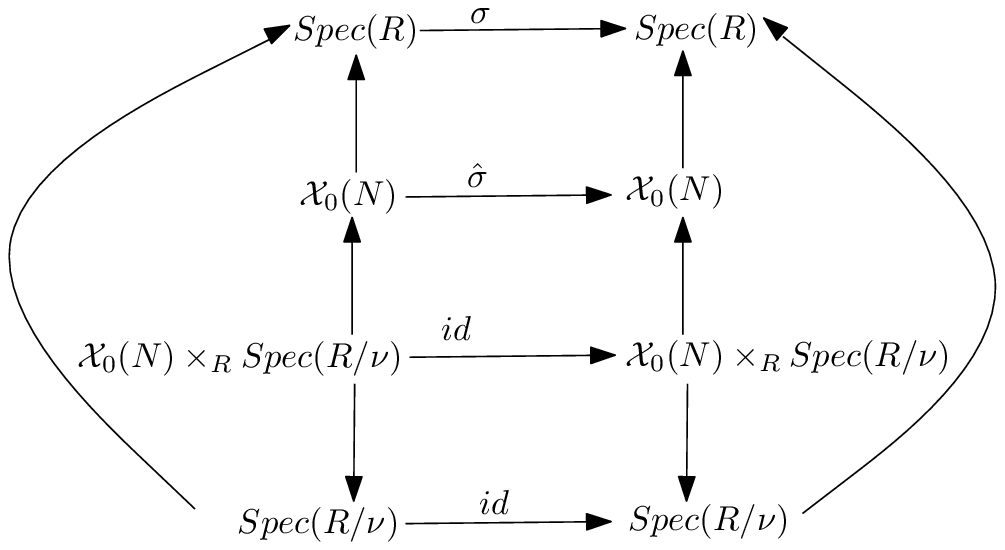}
	\label{fig:ortadiagram}
\end{figure}

We now add to the picture the Atkin-Lehner involution $w_N$  which is a morphism of $Spec(R)$-schemes.

\begin{figure}[h]
	\centering
		\includegraphics{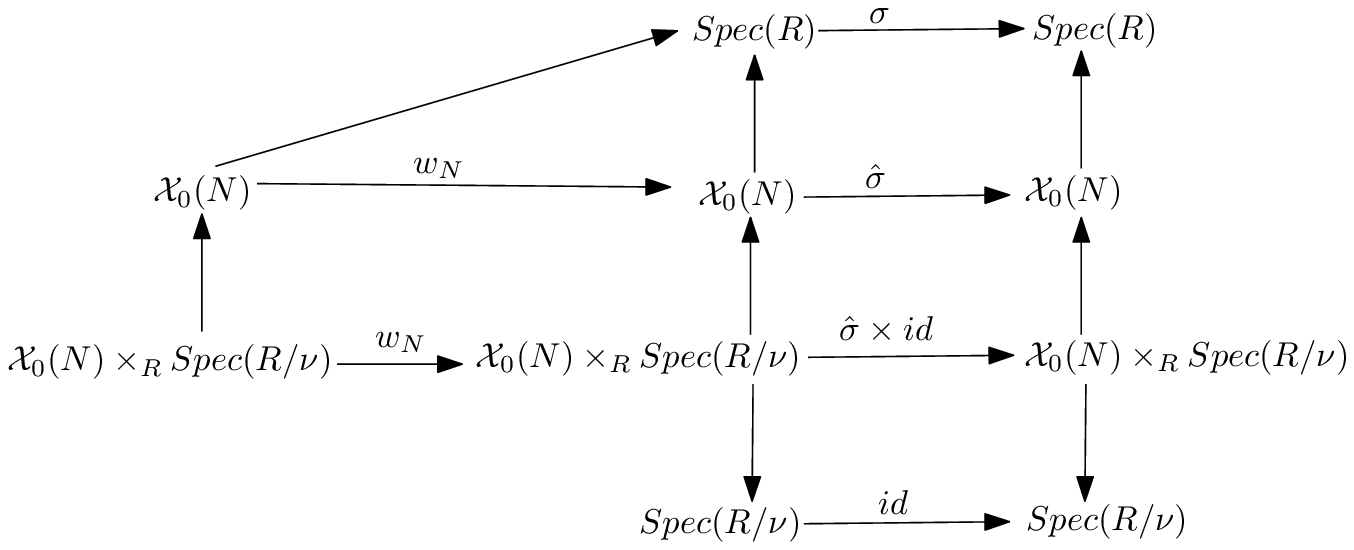}
	\label{fig:buyukdiagram}
	
\end{figure}

  
  
Every point on $X_0(N)(\K_{\nu})$ extends to a morphism $\phi:Spec(R) \rightarrow \mathcal{X}_0(N) $ by properness, and if the point of $X_0(N)(\K_{\nu})$ is fixed by $w_N \circ \sigma$ then the morphism $\phi$ is preserved under composition with $w_N \circ \hat{\sigma}$. To be more precise, let $x$ be a point in $X_0(N)(\K_{\nu})$ such that $w_N\circ \sigma (x)=x$. By properness, $x=\phi \circ i$ where is the injection $i:Spec(\K_{\nu}) \rightarrow Spec(R)$, then $w_N \circ \hat{\sigma} \circ \phi \circ i= \phi \circ i$, hence $w_N \circ \hat{\sigma} \circ \phi= \phi $.

And the diagram shows that the restriction of $\phi$ to the special fiber $\tilde{p}:Spec(R/\nu) \rightarrow \mathcal{X}_0(N) \times_R Spec(R/\nu) $ is a $w_N$-fixed point on $X_0(N)(\F_p)$.
 \end{proof}

 In fact the Proposition \ref{prop:diag} is true even if $p$ is ramified in $\Q(\sqrt{-N})$, in order to have a $K_{\nu}$ rational $w_N \circ \sigma$- fixed point there must be a $w_N$-fixed $\F_p$-rational point of $X_0(N)$. However the converse can not be concluded using Proposition \ref{prop:onlynice}. Since if $p$ is ramified in $\Q(\sqrt{-N})$, $p$ is ramified in $H/\Q$, it is not immediately clear how $B$ ramifies at primes over $p$.
Nonetheless, we have the following result for any prime $p$ ramified in $K$, without any restriction on the decomposition of $p$ in $\Q(\sqrt{-N})$:

\begin{proposition}\label{prop:ramboth}
Let $p$ be a prime ramified in $K$. If $X^d(N)(\Q_p) \neq \emptyset$ then there is a $w_N$-fixed point in $X_0(N)(F_p)$. 
\end{proposition}

It remains to determine when there are $w_N$-fixed points of $X_0(N)(\F_p)$. Let $S_N$ be the set of primes $p$ such that there is a $w_N$ fixed, $\F_p$-rational point on the special fiber of $\mathcal{X}_0(N)_{/R}$. In Proposition \ref{propmain}, we will describe the set $S_N$ explicitly as a Chebotarev set. In addition to the notation introduced in the beginning of the section, $B$ denotes $\Q(j(O))$ where $j(O)$ is $j$-invariant of the order $O=\Z[\sqrt{-N}]$ and $\M:=\Q(\sqrt{-N})$.

 
  
  \begin{proposition}\label{propmain}
  Let $p$ be an odd prime and $\mathcal{P}$ be a prime of $\M$ lying over $p$. Then $p$ is in $S_N$ if and only if There exists a prime $\nu$ of $\B$ lying over $p$ such that $f(\nu|p)$=1 and $\mathcal{P}$ totally splits in $\H/\M$.
  
 
 \end{proposition}
 
 \begin{proof}
 We know that, $\H/\M$ is an abelian extension with Galois group $G$ isomorphic to the ideal class group of $\M$ and $[\H:\M]=[B:\Q]$. The extension $\H/\Q$ is Galois with Galois group $G \rtimes \Z/2\Z$ and the $\Z/2\Z$-fixed subfield of $\H$ is $\B$.  

We want to know for which primes there is a $w_N$-fixed point on $X_0(N)(\F_p)$.  We have shown in Proposition \ref{prop:onlynice} that this is equivalent to the presence of a $w_N$-fixed point on $X_0(N)(\Q_p)$. 

Let $P$ be a $w_N$-fixed point of $X_0(N)$, defined over $B$. Then $P$ reduces to an $\F_p$ point on the special fiber if and only if it is fixed by Frobenius. Recall that since $B/\Q$ is unramified at $p$, Frobenius acts on $B$.  Hence we should find out for which $p$, there exists a prime $\nu$ of $B$ such that $f(\nu|p)=1$. 

Let $\pi_p$ be the Frobenius at $p$. The map $\pi_p$ gives a conjugacy class $\lambda$ in $\Gal(\H/\Q)$ via Artin symbol .


Now we need to find the conjugacy classes of $G \rtimes \Z/2\Z$. Let $g_i \in G$ and $a,b \in \Z/2\Z$, then $(g_1,a)*(g_2,b)=(g_1+\phi_a(g_2),a+b)$ where $\phi_1(g)=-g$ and $\phi_0$ is identity for all $g \in G$. 



Let $(g,a)$ be in $G \rtimes \Z/2\Z$. If $a=0$, $(g,a)$ is only conjugate to itself and $(-g,a)$. If $a=1$ then $(g,a)$ is conjugate to $(-2x+g,a)$ for some $x \in G$.

The conjugacy classes of $G \rtimes \Z/2\Z$ are thus as follows:
\begin{enumerate}
\item $\{(g,0)\}$ one for each $g \in G[2]$.
\item $\{(g,0),(-g,0)\}$ one for each $g$ in $G-G[2]$.
\item $\{(g+2x,1)|x \in G\}$ one for each representative $g$ of $G/2G$.
\end{enumerate}

Let $\pi_p$ be the conjugacy class in $\Gal(\H/\Q)$ given by a Frobenius at $p$. A prime $\nu$ of $\B$ over $p$ has $f(\nu|p)=1$ if and only if the conjugacy class $\pi_p$ contains an element of the form $(0,y)$ for some $y$ in $\Z/2\Z$. Hence, only allowed conjugacy classes are the trivial class and one of the classes of type (3). 

Hence $p$ is in $S_N$ if and only if $\pi_p$ contains an automorphism which fixes $B$, which is equivalent to say that $\mathcal{P}$ totally splits  in $\H/\M$. 


\end{proof}

\begin{remark} 
Note that if $p$ splits in $\M/\Q$ then there are $2$ primes of $\M$ lying over $p$. If $\mathcal{P}$ of $\M$ lying over $p$ splits totally in $\H/\M$ then $p$ splits totally in $\H/\Q$, hence it doesn't matter which prime of $\M$ lying over $p$ we take.
\end{remark}

\begin{remark} Proposition \ref{propmain} determines for which $p$ the field of definition of an elliptic curve whose endomorphism ring contains   $\Z[\sqrt{-N}]$ embeds into $\Q_p$. Then using Proposition \ref{CM}, we get a $w_N$-fixed $\Q_p$-rational point of $X_0(N)$.
\end{remark}

We have thus established a complete criterion for the non-emptiness of $X^d(N)(\Q_p)$, where $p$ is an odd prime ramified in $\K$ but not in $\Q(\sqrt{-N})$.

For $\Q_2$-points, we can argue as follows. Let $d$ and $N$ be square-free integers such that $d \equiv 2,3$ and $-N \equiv 1$ mod $4$. 
 Over $F_2$ there are $2$ elliptic curves one with endomorphism ring $\Z[\frac{1+\sqrt{-7}}{2}]$ ,the ordinary one, and the other supersingular one whose endomorphism ring is Hurwitz quaternions, $B(\Z):=\Z+i\Z+j\Z+\frac{1+i+j+k}{2}\Z$. It is a maximal order in the quaternion algebra ramified only at $2$ and infinity. Hence, if $w_N$-fixed point $(E,C)$ of $X_0(N)(\F_2)$ is ordinary -in particular $N=7$- then $E$ can be lifted to an elliptic curve over a number field $B$ that has complex multiplication by the maximal order of $\Q(\sqrt{-7})$ by Theorem \ref{thm:DeuringLifting}. If $w_N$-fixed point $(E,C)$ of $X_0(N)(\F_2)$ is the supersingular one, then the maximal order of $\Q(\sqrt{-N})$ embeds in $End(E)$ since the local order $B(\Z)\otimes_{\Z}\Z_2$ contains all elements of $B(\Q) \otimes \Q_2$ with norm in $\Z_2$. Therefore by Theorem \ref{thm:DeuringLifting} $E$ can be lifted to an elliptic curve over a number field $B$ that has complex multiplication by the maximal order of $\Q(\sqrt{-N})$. Hence we proved the following lemma:
 
 \begin{lemma}\label{lem:liftingtwo}
 Let $(E,C)$ be a $w_N$-fixed point of $X_0(N)(\F_2)$. Then $E$ can be lifted to an elliptic curve $\tilde{E}$ over a number field $B$ such that $\tilde{E}$ has complex multiplication by the maximal order of $\Q(\sqrt{-N})$.
 \end{lemma}

Let $\tilde{E}$ has CM by the maximal order of $\Q(\sqrt{-N})$. Since $2$ is unramified in $\Q(\sqrt{-N})$ and the Hilbert class field is an unramified extension of $\Q(\sqrt{-N})$, $2$ is unramified in $B/\Q$. This implies that $B_{\nu}$ embeds in $\Q_2$ where $\nu$ is a prime of $B$ lying over $2$. Therefore $\tilde{E}$ induces a point in $X^d(N)(\Q_2)$. Hence we conclude that $X^d(N)(\Q_2)\neq \emptyset$ if and only if there is a $w_N$-fixed point on $X_0(N)(\F_2)$ and such a point exists if and only if $p$ is in the set $S_N$ defined above, exactly as in the case of odd primes. This gives us the following result :


\begin{theorem}\label{thm:ram}
Let $p$ be a prime ramified in $\Q(\sqrt{d})$ and $N$ a squarefree integer such that $p$ is unramified in $\Q(\sqrt{-N})$. Then $X(\Q_p) \neq \emptyset$ if and only if $p$ is in the set $S_N$ defined in Proposition \ref{propmain}.
\end{theorem}


\begin{definition} A set of rational primes $S$ is a \emph{Chebotarev set} if there is some finite normal extension $L/\Q$ such that $p$ is in $S$ if and only if the Artin symbol of $p$ is in some specified conjugacy class or union of conjugacy classes of the Galois group $\Gal(L/\Q)$.  
\end{definition}

The density of a Chebotarev set is well-defined by Chebotarev density theorem, and the set of primes $S_N$ defined in Proposition \ref{propmain} is a Chebotarev set with density $\frac{|2G|+1}{2|G|}$.


\begin{theorem}[Serre, Theorem 2.8 in \cite{serreden}] \label{thm:serreden}
Let $0<\alpha <1$ be  Frobenius density of a set of primes $S$ and $N_S(X):=\{n\leq X: n=\Pi_i p_i, p_i \in S\}$. Then $|N_S(X)| = c_S\frac{X}{\log^{1-\alpha}X}+ O\left (\frac{X}{\log^{2-\alpha}X} \right )$ for some positive constant $c_S$. 
\end{theorem}


\begin{corollary}
Given $N$ let $A$ be the set of squarefree integers $d$ in $[1 \cdots X]$ such that $(d,N) = 1$ and $X^d(N)(\Q_p)$ is nonempty for all primes $p$.
Then $|A|=M_{S_N} \frac{X}{\log^{1-\alpha}X}+ O\left (\frac{X}{\log^{2-\alpha}X} \right )$ where $\alpha=\frac{|2G|+1}{2|G|}$ is the density of $S_N$ and $M_{S_N}$ a positive costant.
\end{corollary}


\begin{proof}
Follows from Theorem \ref{thm:ram} combined with Theorem \ref{thm:serreden}.
\end{proof}

The following corollary gives an answer to a question of Clark asked in \cite{Clark}.

\begin{corollary} \label{cor:Clark}Let $N$ be a prime congruent to $1$ mod $4$, $p$ be an odd prime such that $(N/p)=-1$ and $d=p^*$ then 
\begin{enumerate}
	\item $X^d(N)(\Q_N)=\emptyset$
	\item $X^d(N)(\Q_p)=\emptyset$
	\item $X^d(N)(\Q_\ell) \neq \emptyset$ for any other prime $\ell$ different than $p$ and $N$.
\end{enumerate}

\end{corollary}

\begin{proof}  

The first conclusion was also given in \cite{Clark} and can be seen as a result of Theorem \ref{thm:IDB}. For the second conclusion we proceed as follows:

Let $\M:=\Q(\sqrt{-N})$. Since $N \equiv 1$ mod $4$, the genus field of $\M$ is $\Q(i,\sqrt{N})$. 
Note that since $-N \equiv 3$ mod $4$, ring class field of $\Z[\sqrt{-N}]$ is the Hilbert class field. Let $\B:=\Q(j(\Z[\sqrt{-N}]))$. Since $j(\Z[\sqrt{-N}])$ is  real, $\B \cap \Q(i, \sqrt{N})$ is $\Q$ or $\Q(\sqrt{N})$. Say  $\B \cap \Q(i, \sqrt{N})$ is $\Q$, this implies that the class number of $\Z[\sqrt{-N}]$ is $1$, contradiction.

By Theorem $\ref{thm:ram}$ and Lemma \ref{propmain}, $X^d(N)(\Q_p)=\emptyset$ if and only if $p \notin S_N$. Since $(N/p)=-1$ this is equivalent to say that for all primes $\nu$ of $\B$ lying over $p$, $f(\nu|p)>1$. 

The third conclusion can be derived from Theorem $\ref{thm:inertnotdividing}$.

\end{proof}



   
\section{Examples} \label{sec:exam}
\begin{example}
This example deals with the case $N=7*19$, $p=19$ and $d=3$. According to Theorem \ref{q}, since $(7*19,3)=-1$, $X^d(7*19)(\Q)=\emptyset$. Since $(\frac{3}{19})=-1$ and $19$ divides $N$ we are in the first case of Section \ref{sec:inert}. Using the related theorem in this case(Theorem \ref{thm:IDB}), we can in fact see that $X^d(7*19)(\Q_{19})=\emptyset$. 

To explain this more concretely, we use the idea of the proof Theorem \ref{thm:IDB}. Using the genus formula; $$g(N)=2g(N/p)+n-1$$ where $n$ is the number of supersingular points on $X_0(N)_{\F_p}$, we see that there are 12 supersingular points on $X_0(19*7)_{\F_{19}}$. And there are 2 supersingular j-invariants over $\bar{\F_{19}}$, one of them is 1728(since 19 is congruent 3 mod 4) and the other one is $7$. The elliptic curve $E$ with j-invariant 1728, has 8 order-7 subgroups 
namely $C_1, ..., C_8$. The extra automorphism($[i]$, where $i$ is the fourth root of unity) sends $C_k$ to $C_{9-k}$. Hence we get 4 supersingular points on $X_0(7)_{\F_{19}}$ each with automorphism group isomorphic to $\{\pm 1\}$; namely $P_1=(E,C_1), P_2=(E,C_2), P_3=(E,C_3), P_4=(E,C_4)$. For the elliptic curve with j-invariant $7$, we get 8 many order 7 groups 
each with the smallest possible automorphism group. Hence we get 8 supersingular points, satisfying the total number 12 found above. Since the number of blow-ups is 0, the model of $X_0(7*19)_{\F_{19}}$ is regular and the only $\F_{19}$-rational points are the singular ones, hence $X^d(7*19)(\Q_{19})=\emptyset$, this is also consistent with the observation that $(7*19,3)_{19}=-1$.

Since $(7*19,3)_{19}=-1$, there should be at least one more local Hilbert symbol that is negative, in this case it is $(7*19,3)_{7}=-1$. Now we will do the same operations as above, this time for $X_0(7*19)_{\F_7}$. Using the genus formula, there are 10 supersingular points on $X_0(7*19)_{\F_7}$ and there is one supersingular elliptic curve over characteristic 7, the one with j-invariant $1728$. Let $E$ be the elliptic curve with $j$ invariant $1728$. 
There are 20 tuples like $(E,C_i)$  and each has an automorphism group equal to $\{\pm 1\}$ or primitive fourth roots of unity. Since there are 10 supersingular points, all has the smallest automorphism group i.e. $X^d(7*19)(\Q_7)=\emptyset$ using the same argument above. This also can be deduced direclty from Theorem \ref{thm:IDB}.
\end{example}  

\begin{example} Let $d=5$ and $ N=29$. According to Theorem \ref{q}, since $(5,29)=1$, we have the necessary condition for the existence of a $\Q$-curve of degree $29$ over $\K=\Q(\sqrt{5})$, but it is not guaranteed. Note that this case is not covered by Theorem \ref{thm:Clark}.

If we use Theorem \ref{thm:ram} for the ramified prime $5$, we see that $X^5(\Q_5)$ is empty, hence there is no $\Q$-curve of degree $29$ over $\K$. In order to $X^5(\Q_5)$ to be non-empty, $5$ should split in $H/\Q(\sqrt{-29})$, where $H$ is the Hilbert class field of $\Q(\sqrt{-29})$. However the prime $P|5$ of $\Q(\sqrt{-29})$, decomposes as $P_1P_2$ where inertia degree of $P_i$ is  $3$, hence $X^5(\Q_5)=\emptyset$. Note that $X^5(\Q_p)(29) \neq \emptyset$ for any other prime $p$ different than $5$ using Theorem \ref{thm:ALL}.
\end{example} 

\section{Further Directions}\label{sec:further}

Using Theorem \ref{thm:ALL}, we can produce lots of rational curves which have local points everywhere. One natural follow-up question would be asking about the $\Q$-rational points. There is an answer to this question in the case of imaginary quadratic fields $\K$ and when $N$ is inert in $\K$ implied by the following theorem of Mazur (\cite{Mazur}): 

\begin{theorem}If $\K$ is a quadratic imaginary field, and $N$ is a sufficiently large prime which is inert in $\K$, then $X_0(N)(K)$ is empty. In particular, there are no $\Q$-curves over $\K$ of degree $N$. \label{result}
\end{theorem}   
  
 When $N$ splits in $\K$ and $\K$ is imaginary quadratic or $N$ is inert in $\K$ and $\K$ is real quadratic field, using the formula of Weil given in \cite{Li}, every cuspform associated with a quotient of jacobian of $X^{d}(N)$ has odd functional equation. Thus, conjecturally none of these quotients have Mordell-Weil rank $0$ and we cannot hope to apply Mazur's techniques. The future plan is to prove a result about existence of rational points on $X^d(N)$ where $\K$ is a real quadratic field and $N$ splits in $\K$ using Mazur's techniques.
 
 
 Another direction to go is understanding the reasons of violations to Hasse principle. Say for some $d$ and $N$, $X^d(N)$ has local points for every $p$ but no global points, hence it violates Hasse principle. What is the reason for that?  One natural guess would be Brauer-Manin obstruction. It is a folklore conjecture that for curves, \emph{`Under the assumption that $Sha(J)$ is finite, Brauer-Manin obstruction is the only obstruction to Hasse principle'}, which is a theorem in the following cases: 
\begin{itemize}
\item \emph{(Manin)} If $C$ is proper, smooth of genus $1$ and $Sha(J)$ is finite (\cite{Flynn})
\item \emph{(Scharaschkin)} If $C$ is proper, smooth, has a rational divisor class of degree 1 and $J(\Q)$ and $Sha(J)$ are finite (\cite{sch})
\item \emph{(Skorobogatov)} If $C$ is proper, smooth, $Sha(J)$ is finite and $C$ has no rational divisor class of degree 1 (\cite{sko})

 \end{itemize}
 
 For some specific examples, we will use the following idea of Scharaschkin (\cite{sch}) to see if violation to Hasse principle can be explained by Brauer-Manin obstruction.

Given a smooth, projective, geometrically integral curve $C$ over $\Q$ and a finite set $S$ of good primes, if images of $red$ and $inj$ in the following map do not intersect then $C(\Q)=\emptyset$.

$$\begin{array}{ccc}
C(\Q) & \longrightarrow & Jac(\Q) \\
\downarrow & & \downarrow_{red} \\
\prod_{p \in S} C(\F_p) & \stackrel{inj}{\longrightarrow} & \prod_{p \in S} Jac(\F_p)\\
	
\end{array}$$
 
 In order to apply Scharaschkin's technique, one needs an equation of the curve $C$, generators of $Jac(\Q)$ and also existence of a $\Q$-rational degree one divisor class. In the case of quadratic twists of $X_0(N)$, if the curve is hyperelliptic and $w_N$ is the hyperelliptic involution then finding the equation of the twist is easy. Moreover there exists relatively simple equations of $X_0(N)$ given by Galbraith in \cite{gal} that makes the computations feasible. Such equations for hyperelliptic modular curves is first given by Gonzalez \cite{Gonzalez}, see also the works of Murabayashi \cite{Murabayashi} and Hasegawa \cite{Hasegawa}. It's a result of Ogg (\cite{ogg}) that automorphism group of $X_0(N)$ is $\{1,w_N\}$ for prime $N$ such that $N \neq 2,3,5,7,13,37$. Then the equation of the twist $X^d(N)$ is $dy^2=f_{2g+2}(x)$ where $g$ is the genus of the curve and $f_m$ is a degree $m$ polynomial.
 
 Since for genus $1$ the claim is already proved, we'll restrict to the cases $g \geq 2$ and we want an hyperelliptic curve with $w_N=-1$. Smallest such $N$ is $23$ and $X_0(23)$ is given by $(x^3-x+1)(x^3-8x^2+3x-7)$ in \cite{gal}. The following computations are done using the computer package MAGMA.

 \begin{example} Let $N=23$. We will study the twists of $X_0(23)$ for all primes $d$ between $-300$ and $300$. There are 124 primes between -300 and 300. The twist is given by the equation $y^2=d(x^3-x+1)(x^3-8x^2+3x-7)$. Let $a_1,a_2,a_3$ be roots of $x^3-x+1$, and $P_i=(a_i,0)$ be the corresponding points on $X^d(N)$ then $D=[P_1+P_2+P_3-\infty_1-\infty_2]$ is a rational divisor of degree one on $X^d(N)$.(Similar construction is in \cite{Flynn}).
\begin{enumerate}
\item Let $|d|$ be prime different then $23$ and between $-300$ and $300$, $X^d(N)$ has local points everywhere for 39 values of $d$. Among the 39 values of such $d$, for 10 of them $X^d(N)$ has points with small height(less than 1000) when we eliminate these, we are left with the set: $\{ -283, -271, -263, -251, -227, -223, -211,-199,$

$-191, -83, -59, 17, 37, 53, 61, 89, 97, 101, 109, 113, 137, 149, 157,$ 

$173, 181, 229, 241, 281, 293 \}$

\item The remaining thing we need in order to apply Scharaschkin's technique are the generators of $J^d(\Q)$ where $J^d$ is the Jacobian of $X^d(N)$. This seems to be the hardest thing to do. First lets consider the only one that we were able to apply Scharaschkin's technique.

Say $d=17$: Let $C$ be $X^{17}(23)$ and $J^{17}$ be its jacobian. It can be computed that $J{17}$ has no nontrivial torsion and rank of $J^{17}(\Q)$ is less than or equal to $2$. After a short search we come up with the generators of $J^{17}(\Q)$, given by $D_1=<x^2 + 3, 17*x - 34, 2>$ and 
$D_2=<x^2 - 3/4*x + 5/8, 153/16*x - 17/32, 2>$. The notation means that $D_1=[P_1+\bar{P_1}-\infty_1-\infty_2]$ with $P_1=(a,17*a-34)$ where $a$ is one of the roots of $x^2+3$ and $\bar{P_1}=(\bar{a},17*\bar{a}-34)$. Similarly for $D_2$. 

To apply Scharaschkin's idea, we need a set of primes $S$, in our case $S$ will be $\{13\}$. 

Let $\tilde{C}, \tilde{J^{17}}, \tilde{D_1}, \tilde{D_2}$ be the reductions of $C,J^{17},D_1,D_2$ modulo $13$. Order of $D_1$ and 
$D_2$ are $11$ in $\tilde{J^{17}}(\F_{13})$. 

Say $P$ is in $C(\Q)$ then the reduction $\tilde{P}$ in  
$\tilde{C}(\F_{13})=\{ (1 : 1 : 0), (1 : 12 : 0),(3 : 11 : 1), (3 : 2 : 1), (6 : 8 : 1), (6 : 5 :1), (7 : 3 : 1), (7 : 10 : 1) \}$ and 
$[\tilde{P}+\tilde{P}-\infty_1-\infty_2]$ will be in $\tilde{J^{17}}(\F_{13})$. 

The points in $\tilde{J^{17}}(\F_{13})$ are given by $n_1$ times 
$\tilde{D_1}=[(6,3)+(7,7)-\infty_1-\infty_2]$ and $n_2$ times $\tilde{D_2}=[(a,12a+8)+(b,12b+8)-\infty_1-\infty_2]$ 
($a,b$ are roots of $x^2 - 3/4*x + 5/8$) where $n_1,n_2$ are between $0$ and $10$. If we construct all these linear combinations we see that 
none of these is of the form $[\tilde{P}+\tilde{P}-\infty_1-\infty_2]$ for $\tilde{P}$ in $\tilde{J^{17}(\F_{13})}$, hence there is no $P$ in $C(\Q)$.

\end{enumerate}
\end{example}

This example shows that, the twisted modular curve $X^{17}(23)$ has local points everywhere but no global points, hence it violates Hasse principle and this violation can be explained by Brauer-Manin obstruction.

This example is also interesting in the sense that $23$ is inert in the quadratic field $\Q(\sqrt{17})$. Then, using the formula of Weil given in \cite{Li}, every cuspform associated with a quotient of jacobian of $X^{17}(23)$ has odd functional equation. Thus, conjecturally none of these quotients have Mordell-Weil rank $0$ and Mazur's methods cannot be applied. In fact, the jacobian of $X^{17}(23)$ is simple since $J^{17}(23)$ is an abelian surface with rank $2$ and therefore its only notrivial quotients are the elliptic ones. However, the $q$-expansion of the corresponding newforms of level $23$ has conjugate coefficients in $\Q(\sqrt{-5})$ (see \cite{stein}) hence there is only one isogeny class, therefore $J^{17}(23)$ is simple.

If we continue our search for the same level, $N=23$, we get some more $d$ values such that corresponding twists fails to have rational points.

\begin{example}(continued from above) For $d=173,-211,101,-59,-223$ the rank of $Jac(\Q)$ is $0$ since the 2-Selmer group is trivial. Moreover, the torsion part of $Jac(\Q)$ is also trivial in all these cases. Say there exists $P \in C(\Q)$ where $C$ is the twist $X^d(23)$, then $[P]-D$, where $D$ is as above, is in $J(\Q)$. It can be proved that $D$ is not equivalent to $[P]$ for any point $P$ hence $[P]-D$ is nonzero, contradiction. Therefore $X^d(23)(\Q)=\emptyset$ for $d=173,-211,101,-59,-223$.
\end{example}


\end{document}